\documentclass[journal]{IEEEtran}
\usepackage{amsmath,amssymb}
\usepackage{graphicx}
\usepackage{overpic}
\usepackage{tikz}
\usepackage{bm}

\newtheorem{example}{Example}
\newtheorem{algorithm}{Algorithm}

\newcommand{\K}{\mathcal K}
\newcommand{\M}{\mathcal M}
\newcommand{\mcO}{\mathcal O}
\newcommand{\U}{\mathcal U}
\newcommand{\bA}{\bm A}
\newcommand{\bB}{\bm B}

\newcommand{\bg}{\bm g}
\newcommand{\bG}{\bm G}
\newcommand{\bh}{\bm h}

\newcommand{\bu}{\bm u}
\newcommand{\bU}{\bm U}
\newcommand{\bv}{\bm v}
\newcommand{\bV}{\bm V}
\newcommand{\bw}{\bm w}
\newcommand{\bW}{\bm W}

\newcommand{\bX}{\bm X}
\newcommand{\by}{\bm y}
\newcommand{\bY}{\bm Y}
\newcommand{\bz}{\bm z}
\newcommand{\bZ}{\bm Z}
\newcommand{\bgamma}{\bm\gamma}
\newcommand{\bGamma}{\bm\Gamma}
\newcommand{\bphi}{\bm\phi}
\newcommand{\bPhi}{\bm\Phi}
\newcommand{\bSigma}{\bm\Sigma}
\newcommand{\bLambda}{\bm\Lambda}
\newcommand{\bOmega}{\bm\Omega}
\newcommand{\rar}{\rightarrow}
\newcommand{\C}{\mathbb C}
\newcommand{\R}{\mathbb R}

\begin{document}
\title{Dynamic mode decomposition for interconnected control systems}
\author{Byron~Heersink,
		Michael~A.~Warren,
		and Heiko~Hoffmann
\thanks{This work was funded under DARPA contract N66001-16-C-4053. The views expressed are those of the authors and do not reflect the official policy or position of the Department of Defense or the U.S. Government. Distribution Statement `A': Approved for Public Release, Distribution Unlimited.}
\thanks{B.\ Heersink was with HRL Laboratories, LLC, Malibu, CA 90265. He is now with the Department of Mathematics, The Ohio State University, 231 W.\ 18th Ave., Columbus, OH 43210 (e-mail: heersink.5@osu.edu).}
\thanks{M.\ A.\ Warren and H.\ Hoffmann are with HRL Laboratories, LLC, 3011 Malibu Canyon Rd., Malibu, CA 90265 (email: mawarren@hrl.com; hhoffmann@hrl.com).}}

\maketitle

\begin{abstract}
Dynamic mode decomposition (DMD) is a data-driven technique used for capturing the dynamics
of complex systems. DMD has been connected to spectral analysis of the Koopman operator, and
essentially extracts spatial-temporal modes of the dynamics from an estimate of the Koopman
operator obtained from data. Recent work of Proctor, Brunton, and Kutz has extended DMD and
Koopman theory to accommodate systems with control inputs: dynamic mode decomposition
with control (DMDc) and Koopman with inputs and control (KIC). In this paper, we introduce a
technique, called Network dynamic mode decomposition with control, or Network DMDc, which
extends the DMDc to interconnected, or networked, control systems. Additionally, we provide an
adaptation of Koopman theory for networks as a context in which to perform this algorithm. The
Network DMDc method carefully analyzes the dynamical relationships only between components
in systems which are connected in the network structure. By focusing on these direct dynamical
connections and cutting out computation for relationships between unconnected components,
this process allows for improvements in computational intensity and accuracy.
\end{abstract}

\begin{IEEEkeywords}
Koopman operator theory; System identification; Network analysis and control; Computational methods.
\end{IEEEkeywords}

\section{Introduction}\label{sec:intro}

\IEEEPARstart{D}{ynamic} mode decomposition (DMD) is a method developed by
Schmid and Sesterhenn \cite{SS,S1} used in the model reduction and decomposition
of complex dynamical systems. This data-driven method is performed on time series data
of a given system and attempts to identify a linear model for the dynamics, which is ideally
of reduced order. Then the prominant behavior of the system is extracted from the linear
model's eigenvectors, or ``modes'', whose dynamics are governed simply by their
corresponding eigenvalues. 

It has been shown that DMD is strongly related to Koopman operator theory.
Originally defined by Koopman in 1931 \cite{HSTHS}, the Koopman operator
is a linear infinite-dimensional operator on the space of observables of any
dynamical system, including nonlinear systems. The work of Mezi\'c \cite{M1}
was the first to apply the spectral analysis of the Koopman operator to the
model reduction of systems. Later, Rowley et al.\ \cite{RMBSH} fundamentally
linked DMD to Koopman theory by showing that DMD is essentially a method
of approximating the Koopman operator and its spectral decomposition.

In recent years, DMD and Koopman theory has been applied to a variety of
different fields. It's first and most notable application area is fluid dynamics.
See, e.g., \cite{RMBSH,S1,S2,TRLBK}, as well as \cite{M2} for a review.
Other applications include power systems \cite{SM1,SM2,SM3,SMRH},
video processing \cite{GK,KGB,ED}, epidemiology \cite{PE}, robotics \cite{BSVJA},
neuroscience \cite{BJOK}, and finance \cite{MK}. Additionally, there has been
a lot of effort to refine and advance DMD and Koopman theory
themselves. In particular, DMD has been improved by Tu et al.\ \cite{TRLBK}
to what is recognized as its preeminant form. Examples of innovations building
upon DMD and Koopman theory include extended and kernel DMD \cite{WKR,WRK},
which seek to more accurately approximate the Koopman operator by incorporating
measurements of appropriate nonlinear observables explicitly or implicitly through
use of a kernel; multi-resolution DMD \cite{KFB}; and the incorporation of sparsity
\cite{JSN}, compression \cite{BPK}, and de-biasing \cite{HRDC}.

Another recent extension of DMD that can be applied to dynamical systems with
control inputs is the dynamic mode decomposition with control (DMDc) \cite{PBK1}.
Similar to DMD, DMDc seeks to find a linear model which approximates the dynamics
of a given control system using only data measurements of the state and inputs of
the system. DMDc has been applied to modeling a rapidly pitching airfoil \cite{DSRW}, and
has been leveraged to produce a generalization of Koopman theory to incorporate
control inputs \cite{PBK2}. Other efforts to extend Koopman theory to control systems
are \cite{BBPK,WHDKR,KM}.

In this paper, we introduce a further extension of DMDc to interconnected, or networked,
control systems. That is, we formulate a DMDc algorithm specialized for analyzing systems
which are composed of smaller subsystems arranged in a network structure, where each
subsystem corresponds to a node in the network, and the edges in the network represent
the dynamical interactions between subsystems. We call this algorithm Network Dynamic
Mode Decomposition with Control, or Network DMDc. We also adapt Koopman theory to
networked control systems to provide a framework within which to apply the Network
DMDc algorithm.

Examples of systems with a network structure include chemical reaction networks,
epidemiological networks capturing the transmission of diseases through different spatial
locations or groups, power system networks, and gene regulatory networks. The main
idea behind the Network DMDc algorithm is to exploit, when possible, existing network
structure in complex systems to yield improvements in computation intensity and
precision over standard DMDc.

In Section \ref{sec:background}, we review background on standard Koopman theory, the
dynamic mode decomposition, and their analogues incorporating inputs and control.
Then in Section \ref{sec:extnetworks}, we extend dynamic mode decomposition and
Koopman theory to networked systems. We then present some examples demonstrating
the network DMDc algorithm and some of its benefits in Section \ref{sec:examples}. Finally,
we provide concluding remarks in Section \ref{sec:conclusion}.

\section{Koopman theory and dynamic mode decomposition with control}\label{sec:background}

\subsection{Koopman theory}\label{ssec:Koopman}

In this section, we outline the
basics of Koopman theory. The reader is referred to
\cite{RMBSH,BMM} for more details.
Let $\M$ be a state space of a discrete dynamical system
$T:\M\rar\M$, whose evolving trajectories are sequences
$x_1,x_2,x_3,\ldots$ in $\M$ such that
\begin{equation}\label{DS}
x_{k+1}=T(x_k).
\end{equation}
We define the Koopman operator as the operator $\K$ that acts on
scalar-valued observable functions $g:\M\rar\R$ according to
\[(\K g)(x)=g(T(x)).\]
We think of $\K$ as a linear operator on a vector space of observables on
$\M$ that we denote by $\mcO(\M)$, and which is commonly chosen
to be a Hilbert space (e.g., the functions on $\M$ which are
square-integrable with respect to some measure). The main idea of the
Koopman method is to analyze the dynamical system \eqref{DS} through
the behavior of the observables on the state space. In particular, the goal is
to use data recording the value of various observables to find the eigenvalues
and eigenfunctions of $\K$, which would then help us better understand the
dynamics of the system. Assume $\varphi_j:\M\rar\R$, $j=1,2,\ldots$, are the
eigenfunctions of $\K$ with corresponding eigenvalues $\lambda_j\in\C$, $j=1,2,\ldots$,
so that
\[(\K\varphi_j)(x)=\lambda_j\varphi_j(x).\]
Then for a given (vertical) vector valued observable $\bg:\M\rar\R^n$, each of
whose components lie in the span of the eigenfunctions, we can write
\[\bg(x)=\sum_{j=1}^\infty\varphi_j(x)\bv_j.\]
 The operator $\K$ then acts on $\bg$ according to
\[(\K\bg)(x)=\sum_{j=1}^\infty\lambda_j\varphi_j(x)\bv_j.\]
The vectors $\bv_j$, which are called the Koopman modes associated
to $\bg$, are the components of $\bg$ whose dynamics can be simply
discerned from the corresponding eigenvalues $\lambda_j$, which contain
growth or decay rates and the oscillation frequency of the modes.

Koopman theory is also defined for continuous-time dynamical systems.
In this context, $\M$ is a subset of Euclidean space and the trajectories
$x=x(t)$ of the system are governed by a differential equation
\[\frac{dx}{dt}=F(x);\]
or more generally, $\M$ can be a smooth manifold and the dynamics
governed by a vector field on $\M$. (In both cases, we assume the
trajectories $x(t)$ are defined for all $t\geq0$.) We then have the flow function
$\Phi:\M\times[0,\infty)\rar\M$ mapping a pair $(x,t)\in\M\times[0,\infty)$ to the
point in $\M$ obtained by following the dynamics of the system for
time $t$ starting at the point $x$. Then a semigroup of operators
$\{U^t:t\geq0\}$ can be defined on observables $f:\M\rar\R$ by
\[(U^tf)(x)=f(\Phi(x,t)).\]
Since DMD-based methods are performed on data at a discrete set of
times, it is useful to view the system as a discrete-time system by
fixing a time increment $\Delta t>0$ and defining the map
$T:\M\rar\M$ by
\[T(x)=\Phi(x,\Delta t);\]
the corresponding Koopman operator $\K$ is then equal to $U^{\Delta t}$.

\subsection{Dynamic mode decomposition}\label{ssec:dmd}

We now describe how the spectral analysis of the Koopman operator can be done
via dynamic mode decomposition. See \cite{TRLBK} for more details.

We consider two sets of (vertical) data vectors
$\{\bz_1,\ldots,\bz_m\},\{\by_1,\ldots,\by_m\}\subseteq\R^n$.
We think of these vectors as measurements on a dynamical system
such that for every $k$, $\by_k$ is the measurement of the system
that follows the measurement $\bz_k$ after a fixed time increment independent
of $k$. Our goal is to try to find a linear model for the dynamics so that
\[\by_k\approx\bA\bz_k\]
for some matrix $\bA$. In other words, defining the matrices
\[\bZ=\big[\bz_1\enspace\bz_2\enspace\cdots\enspace\bz_m\big]\enspace\text{and}\enspace
\bY=\big[\by_1\enspace\by_2\enspace\cdots\enspace\by_m\big],\]
we wish to have
\[\bY\approx\bA\bZ.\]

We define $\bA$ to be
\[\bA=\bY\bZ^\dagger,\]
where $\bZ^\dagger$ denotes the Moore-Penrose pseudoinverse of $\bZ$.
The matrix $\bY\bZ^\dagger$ is an ideal candidate for $\bA$ since it is the
matrix minimizing $\|\bA\bZ-\bY\|_F$, where $\|\cdot\|_F$ denotes the Frobenius
norm. If $\bY=\bA\bZ$, then $\bX=\bA$ is the solution of $\bY=\bX\bZ$
minimizing $\|\bX\|_F$. The dynamic mode decomposition of the pair $(\bZ,\bY)$
is the eigendecomposition of the matrix $\bA$. However, the necessary computations
could be intensive if the system is sufficiently large, in which case one can consider
a reduced-order model for $\bA$. This is done by the following:
\begin{algorithm}[DMD \cite{TRLBK}]
\begin{enumerate}
\item Compute the reduced and appropriately truncated SVD of $\bZ$:
\[\bZ\approx\bU\bSigma\bV^*.\]
\item Define the reduced-model $\tilde\bA$ for $\bA$ by
\[\tilde\bA=\bU^*\bA\bU=\bU^*\bY\bV\bSigma^{-1}.\]
\item Compute the eigendecomposition of $\tilde\bA$:
\[\tilde\bA\bW=\bW\bLambda,\qquad\bLambda=\mathrm{diag}(\lambda_1,\ldots,\lambda_r).\]
\item For each column $\bw$ of $\bW$ with corresponding eigenvalue
$\lambda\neq0$, compute the associated eigenvector $\bphi$ of
$\bA$ according to
\[\bphi=\lambda^{-1}\bY\bV\bSigma^{-1}\bw\]
These $\bphi$ make up part of the eigendecomposition
of $\bA$, forming the columns of a matrix $\bPhi$ such that
\[\bA\bPhi\approx\bPhi\bLambda.\]
\end{enumerate}
\end{algorithm}

In the context of Koopman theory as outlined in the previous section, the
DMD algorithm is applied to the data
\[\bz_k=\bg(x_k),\enspace\by_k=\bg(w_k),\qquad k=1,\ldots,m\]
where $x_1,\ldots,x_m,w_1,\ldots,w_m\in\M$ is such that $w_k=T(x_k)$,
and $\bg$ is a vector of observables as defined above. The
eigendecomposition $(\bLambda,\bPhi)$ resulting from the
DMD then gives an approximation of the eigenvalues and Koopman modes of
the Koopman operator $\K$ with respect to the observable $\bg$.

\subsection{Koopman theory and dynamic mode decomposition with control}\label{ssec:dmdc}

Next, we outline the work of Proctor et al.\ \cite{PBK1,PBK2} in generalizing
Koopman theory and the DMD to incorporate control inputs. Because of the generality
and concrete framework of Koopman theory, we present the dynamic mode
decomposition with control \cite{PBK1} as encompassed within
Koopman theory with inputs and control (KIC) \cite{PBK2}.

We consider a discrete control system $T:\M\times\U\rar\M$, where $\M$ is the state
space of the system as above, and $\U$ is the space of controls. The trajectories
of this system are sequences $x_1,x_2,\ldots$ in $\M$, with corresponding input
sequences $u_1,u_2,\ldots$ in $\U$ such that
\[x_{k+1}=T(x_k,u_k).\]
In \cite{PBK2}, the Koopman operator $\K$ in this context acts on
observables $g:\M\times\U\rar\R$ according to
\[(\K g)(x,u)=g(T(x,u),*),\]
where $*$ can be chosen in different ways depending on how one wishes to treat
the inputs. One can then attempt to analyze $\K$ as an operator on $\M\times\U$
via DMD. Here, however, we restrict the domain of $\K$ to observables
on $\M$ so that for a given observable $g:\M\rar\R$, we have
\[(\K g)(x,u)=g(T(x,u));\]
and we avoid the ambiguity of having to choose $*$ above. We think of $\K$ as a
linear operator mapping a vector space of observables $\mcO(\M)$ on $\M$ to
another vector space of observables $\mcO(\M\times\U)$ on $\M\times\U$.

As in the autonomous case, one can adapt Koopman theory to continuous-time
control systems. In this context, $\M$ and $\U$ are Euclidean spaces and the
dynamics of the state $x=x(t)$ is governed by a differential equation
\[\frac{dx}{dt}=F(x,u),\]
where $u=u(t)$ is an input signal. (One can more generally consider when
$\M$ and $\U$ are smooth manifolds.) Then for each fixed input $u\in\U$,
one has a function $\Phi_u:\M\times[0,\infty)\rar\M$ yielding the flow of
points in $\M$ assuming the input signal is fixed at $u$. These functions in
turn induce operators $\{U_u^t:t\geq0,u\in\U\}$ defined by
\[(U_u^tg)(x)=g(\Phi_u(x,t)).\]
One can then define the Koopman operator $\K$ by
\[(\K g)(x,u)=(U_u^{\Delta t}g)(x)\]
for some fixed $\Delta t>0$. Note that in doing this, we must restrict the
control signals we consider to those which are constant on time intervals
of length $\Delta t$.

Now let $\bg:\M\rar\R^n$ be a vector of observables on $\M$ and $\bh:\U\rar\R^l$
a vector of observables on $\U$. Also let $(x_1,u_1,w_1),(x_2,u_2,w_2),\ldots,(x_m,u_m,w_m)$
be triples such that $w_k=T(x_k,u_k)$ and define the matrices
\begin{align*}
\bZ&=\big[\bz_1\enspace\bz_2\enspace\cdots\enspace\bz_m\big]
=\big[\bg(x_1)\enspace\bg(x_2)\enspace\cdots\enspace\bg(x_m)\big]\\
\bY&=\big[\by_1\enspace\by_1\enspace\cdots\enspace\by_m\big]
=\big[\bg(w_1)\enspace\bg(w_2)\enspace\cdots\enspace\bg(w_m)\big]\\
\bGamma&=\big[\bgamma_1\enspace\bgamma_1\enspace\cdots\enspace\bgamma_m\big]
=\big[\bh(u_1)\enspace\bh(u_2)\enspace\cdots\enspace\bh(u_m)\big]
\end{align*}
The DMDc algorithm can then be applied in an effort to find a matrices $\bA$ and $\bB$ such that
\[\bg(w_k)\approx\bA\bg(x_k)+\bB\bh(u_k).\]
In other words, if we let $\bOmega=\left[\begin{smallmatrix}\bZ\\\bGamma\end{smallmatrix}\right]$, we wish for $\bG=[\bA\enspace\bB]$ to satisfy
\[\bY\approx\bG\bOmega.\]
We define $\bG$ to be
\[\bG=\big[\bA\enspace\bB\big]=\bY\bOmega^\dagger.\]
The DMDc of the triple $(\bZ,\bY,\bGamma)$ is the eigendecomposition
of the matrix $\bA$. For ease of computation
for large systems, the DMDc algorithm can also
compute a reduced-order model for $\bG$, and then the eigendecomposition for
the corresponding reduced-order model of $\bA$, which gives approximate dynamic
modes for the system. This algorithm is described as follows:

\begin{algorithm}[DMDc \cite{PBK1}]
\begin{enumerate}
\item Compute the reduced and appropriately truncated SVD of $\bOmega$:
\[\bOmega\approx\bU\bSigma\bV^*.\]
and let $p$ be the truncation value so that $\bU\in\R^{(n+l)\times p}$,
$\bSigma\in\R^{p\times p}$, and $\bV^*\in\R^{p\times m}$. Note that
\begin{align*}
\bG&\approx\bY\bV\bSigma^{-1}\bU^*\text{, i.e.,}\\
\big[\bA\enspace\bB\big]&\approx\big[\bY\bV\bSigma^{-1}\bU_1^*\quad\bY\bV\bSigma^{-1}\bU_2^*\big],
\end{align*}
where $\bU_1^*\in\R^{p\times n}$ and $\bU_2^*\in\R^{p\times l}$
such that $\bU^*=[\bU_1^*\enspace\bU_2^*]$.
\item Compute the reduced and appropriately truncated SVD of $\bY$:
\[\bY\approx\hat\bU\hat\bSigma\hat\bV^*.\]
and let $r$ be the truncation value so that $\hat\bU\in\R^{n\times r}$,
$\hat\bSigma\in\R^{r\times r}$, and $\hat\bV^*\in\R^{r\times m}$.
\item Compute the reduced-order model $[\tilde\bA\enspace\tilde\bB]$
of $[\bA\enspace\bB]$ as
\begin{align*}
\big[\tilde\bA\enspace\tilde\bB\big]=\big[\hat\bU^*\bY\bV\bSigma^{-1}\bU_1^*\hat\bU\quad\hat\bU^*\bY\bV\bSigma^{-1}\bU_2^*\big].
\end{align*}
(Intuitively, we have $\tilde\bA\approx\hat\bU^*\bA\hat\bU$ and
$\tilde\bB\approx\hat\bU^*\bB$.)
\item Compute the eigendecomposition for $\tilde\bA$:
\[\tilde\bA\bW=\bW\bLambda\]
\item Compute the approximate eigenvectors of $\bA$ associated to
the columns of $\bW$, forming the columns of a matrix $\bPhi$:
\[\bPhi=\bY\bV\bSigma^{-1}\bU_1^*\hat\bU\bW.\]
\end{enumerate}
\end{algorithm}
The reduced-order model for the state space measurements is then
$\R^r$, and we think of the reduced-order measurement corresponding
to $\bz\in\R^n$ as $\tilde{\bz}=\hat{\bU}^*\bz\in\R^r$. On the other
hand, the reduced-order measurement $\tilde{\bz}$ relates to the
original measurement as $\bz=\hat{\bU}\tilde{\bz}$. (Note $\hat{\bU}^*$
is a left inverse for $\hat{\bU}$.)

\section{Extension to networked systems}\label{sec:extnetworks}

Now that we have outlined the DMDc algorithm and Koopman theory with
control inputs, we can explain their extension to networked systems. As 
before, we present Network DMDc in the context of Koopman theory
due to its generality. 

\subsection{Networked control systems}\label{ssec:networksystems}

We begin by defining precisely what we mean by a networked control
system.  Let $G$ be a directed graph with vertices partitioned into two
disjoint sets $N=\{v_1,v_2,\ldots,v_\nu\}$ and
$I=\{e_1,\ldots,e_\mu\}$ such that the only edges connected to
the vertices in $I$ are directed outward from those vertices. We associate
each vertex $w$ in $N\cup I$ with a set $P_w$ that represents a
component of the system. If $w\in N$, then $P_w$ is a component
of the state space of the system, while if $w\in I$, then $P_w$ is a
component of the input space. The entire state space is therefore
$\M:=\prod_{v\in N}P_v$, and the entire input space is
$\U:=\prod_{e\in I}P_e$. (The set $I$ is allowed to be empty, in which
case the system is autonomous.)

For $v\in N$, let $I_v\subseteq N\cup I$ be the set of vertices having
an outgoing edge pointing into $v$. We then have the transition
function $T_v:P_v\times\prod_{w\in I_v}P_w\rar P_v$ which governs
the behavior of component $P_v$: if $(x_w)_{w\in N\cup I}$ are the
state and input components of the system at a particular time, the
$P_v$ component of the state at the following time is
$T_v(x_v,(x_w)_{w\in I_v})$. Thus the edges of $G$ represent the
pattern of influence the different components of the system have on
one another. The individual transition functions can then be composed
to produce the transition function $T:\M\times\U\rar\M$ of the entire
system. Thus the graph $G$, the state and input spaces $P_w$, and
the transition functions $T_v$ define a networked control system.

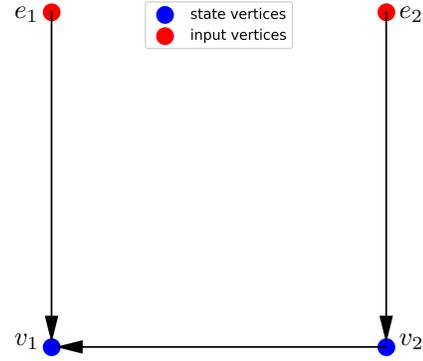
\begin{figure}
\centering
\begin{overpic}[trim=85pt 40pt 85pt 40pt, clip, width=0.3\textwidth]{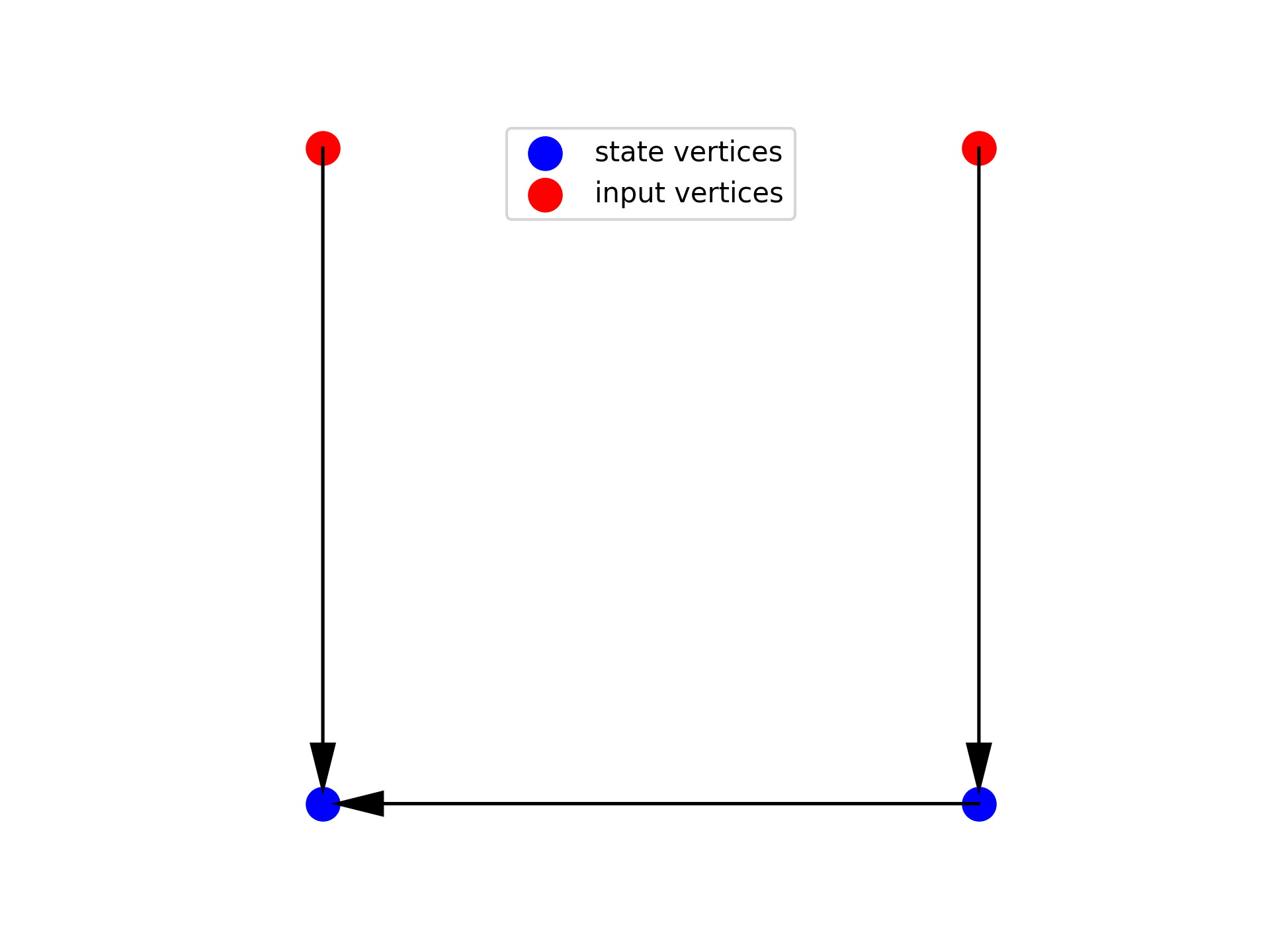}
\put(2,5){$v_1$}
\put(96,5){$v_2$}
\put(2,85){$e_1$}
\put(96,85){$e_2$}
\end{overpic}
\caption{A simple network with $2$ state vertices and $2$ input vertices}
\label{fig:simplenetwork}
\end{figure}

For example, consider the simple network depicted in Figure \ref{fig:simplenetwork}.
It has state vertices $N=\{v_1,v_2\}$ and input vertices $I=\{e_1,e_2\}$.
Since there are edges from $v_2$ and $e_1$ to $v_1$, the transition
function governing the dynamics of the $v_1$-component is of the
form $T_{v_1}:P_{v_1}\times P_{v_2}\times P_{e_1}\rar P_{v_1}$.
Similarly, the edge from $e_2$ to $v_2$ implies that the transition
function for $v_2$ is of the form
$T_{v_2}:P_{v_2}\times P_{e_2}\rar P_{v_2}$. Thus the complete
transition function of the system
$T:\prod_{j=1,2}P_{v_j}\times\prod_{j=1,2}P_{e_j}\rar\prod_{j=1,2}P_{v_j}$
is given by
\[T(x_{v_1},x_{v_2},x_{e_1},x_{e_2})=(T_{v_1}(x_{v_1},x_{v_2},x_{e_1}),T_{v_2}(x_{v_2},x_{e_2})).\]

One can consider analogous continuous-time networked systems as well. Then
all $P_w$ are Euclidean spaces or smooth manifolds, and each transition function
$T_v$ is replaced by a differential equation of the form
\[\frac{dx_v}{dt}=F_v(x_v,(x_w)_{w\in I_v}),\]
or an appropriate map from $P_v\times\prod_{w\in I_v}P_w$ to the tangent bundle
of $P_v$. Such a system can then be approximated by a discrete-time networked
system by a process similar to that given in Section \ref{ssec:dmdc}. Only an
approximation can be made in general due to possible dependence of state
components on other state components that change continuously, which should be
constant in time intervals over which the system is discretized.

\subsection{Koopman theory for networks}\label{ssec:Koopmannets}

In this section, we provide a Koopman theory framework for
the network systems we defined in the previous section. First, for each
$v\in N$, we consider the transition function
$T_v:P_v\times\prod_{w\in I_v}P_w\rar P_v$ as a discrete control system
with state space $P_v$ and input space $\prod_{w\in I_v}P_w$. We can
then define the corresponding Koopman operator $\K_v$ as in Section
\ref{ssec:dmdc} so that for $g:P_v\rar\R$, we have
\[(\K_vg)(x_v,(x_w)_{w\in I_v})=g(T_v(x_v,(x_w)_{w\in I_v})).\]
So we think of $\K_v$ as an operator between vector spaces of observables
$\mcO(P_v)\rar\mcO(P_v\times\prod_{w\in I_v}P_w)$.

Analogous to how the transition functions $T_v$ compose to form $T$,
we wish to see how we can compose the operators $\K_v$ to obtain the
Koopman operator corresponding to $T$. First, by composing it with the
natural inclusion
\begin{equation}\label{oinclusion}
\mcO\left(P_v\times\prod_{w\in I_v}P_w\right)\rar\mcO(\M\times\U),
\end{equation}
(assuming the former is a subspace of the latter) we may view $\K_v$
as an operator from $\mcO(P_v)$ to $\mcO(\M\times\U)$. Then taking
the tensor product of $\K_v$ over $v\in N$ yields an operator
\[\underset{v\in N}{\otimes}\K_v:\bigotimes_{v\in N}\mcO(P_v)\rar\mcO(\M\times\U).\]
We wish to identify the tensor product $\bigotimes_{v\in N}\mcO(P_v)$
with an appropriate space of observables on $\prod_{v\in N}P_v$, ideally
$\mcO(\M)$, via the identification of
$\underset{v\in N}{\otimes}f_v\in\bigotimes_{v\in N}\mcO(P_v)$
with $\prod_{v\in N}f_v$ as a function on $\M$. If this identification can
be made, we can then identify $\underset{v\in N}{\otimes}\K_v$ with
the Koopman operator $\K$ of the entire system.

Some assumptions under which this process works is when for each
$w\in N\cup I$, $\mcO(P_w)=L^2(P_w,\mu_w)$, where $\mu_w$
is a finite measure on $P_w$, and $\K_v$ is bounded. In this case,
the inclusions \eqref{oinclusion} are valid, and the appropriate
completion of $\bigotimes_{v\in N}L^2(P_v)$ is identified with
$L^2(\M,\underset{v\in N}{\otimes}\mu_v)$. (See for instance
\cite[Example 2.6.11]{KR}.)

\subsection{Network DMDc}\label{ssec:netdmdc}

The previous section provides a way of decomposing the Koopman
operator of a networked system into smaller operators. This gives
us a framework in which to formulate our Network DMDc algorithm.
The basic idea is to apply the original DMDc algorithm to each
operator $\K_v$ to get a ``local'' analysis of the behavior of the system
at each state space component $P_v$ in response to the components
it is immediately influenced by. The resulting linear control systems are
then composed to obtain a linear system approximating the whole network.

First of all, let $(x_1,u_1,w_1),\ldots,(x_m,u_m,w_m)$ be triples in
$\M\times\U\times\M$ such that $w_k=T(x_k,u_k)$, and denote the
$P_{v_j}$ component of $x_k$ and $y_k$ as $x_{k,j}$ and $y_{k,j}$,
respectively, and the $P_{e_j}$ component of $u_k$ as $u_{k,j}$. Then
let $\bg:\prod_{j=1}^\nu P_{v_j}\rar\R^n$ and
$\bh:\prod_{j=1}^\mu P_{e_j}\rar\R^l$ be observables of the form
\begin{align*}
&\bg=\left[\begin{array}{c}\bg_1\\\vdots\\\bg_\nu\end{array}\right]\enspace\text{and}\enspace\bh=\left[\begin{array}{c}\bh_1\\\vdots\\\bh_\mu\end{array}\right],\enspace\text{where}\\
&\bg_j:P_{v_j}\rar\R^{n_j},\quad j=1,\ldots,\nu\\
&\bh_j:P_{e_j}\rar\R^{l_j},\quad j=1,\ldots,\mu.
\end{align*}
The first step of the Network DMDc process is to form subsystems centered
at each state vertex. Specifically, for each $v_j\in N$, we consider the
``local subsystem'' consisting of $v_j$, which we think of as the state vertex
of the subsystem, and
\[I_j=\{v_{k_j(1)},v_{k_j(2)},\ldots,v_{k_j(\alpha_j)},e_{\ell_j(1)},e_{\ell_j(2)},\ldots,e_{\ell_j(\beta_j)}\},\]
which we think of as the set of input vertices of the subsystem. Next,
define the observable $\bh^{(j)}$ on $\prod_{w\in I_j}P_w$ by
\[\bh^{(j)}=\left[\begin{array}{c}\bg_{k_j(1)}\\\ldots\\\bg_{k_j(\alpha_j)}\\\bh_{\ell_j(1)}\\\vdots\\\bh_{\ell_j(\beta_j)}\end{array}\right].\]
We then define the triple $(\bZ_j,\bGamma_j,\bY_j)$ by
\begin{align*}
\bZ_j&=\big[\bg_j(x_{1,j})\enspace\bg_j(x_{2,j})\enspace\cdots\enspace\bg_j(x_{m,j})\big]\\
\bY_j&=\big[\bg_j(w_{1,j})\enspace\bg_j(w_{2,j})\enspace\cdots\enspace\bg_j(w_{m,j})\big]\\
\bGamma_j&=\big[\bh^{(j)}(\bu_1^{(j)})\enspace\bh^{(j)}(\bu_2^{(j)})\enspace\cdots\enspace\bh^{(j)}(\bu_m^{(j)})\big],
\end{align*}
where $\bu_k^{(j)}=(x_{k,k_j(1)},\ldots,x_{k,k_j(\alpha_j)},u_{k,\ell_j(1)},\ldots,u_{k,\ell_j(\beta_j)})$.

We analyze the subsystem centered at $v_j$ by applying the DMDc to
$(\bZ_j,\bGamma_j,\bY_j)$. In the case where a reduced-order model
is not sought, the DMDc yields matrices $\bA_{j,j},\bB_j$ that model the
behavior of $\bg_j$ by linear control:
\[\bg_j(x_{k+1,j})\approx\bA_{j,j}\bg_j(x_{k,j})+\bB_j\bh^{(j)}(u_k^{(j)}),\]
which we can rewrite as
\begin{align*}
\bg_j(x_{k+1,j})\approx\text{}&\bA_{j,j}\bg_j(x_{k,j})+\sum_{i=1}^{\alpha_j}\bA_{j,k_j(i)}\bg_{k_j(i)}(x_{k,k_j(i)})\\
&+\sum_{i=1}^{\beta_j}\bB_{j,\ell_j(i)}\bh_{\ell_j(i)}(u_{k,\ell_j(i)}).
\end{align*}
Repeating this process for every $j=1,\ldots,\nu$ gives us a linear control
approximation for each local subsystem, which we can combine to obtain
an approximation for the whole system. Indeed, we can define $A_{j,i}$
to be the zero matrix in $\R^{n_j\times n_i}$ if there is not an edge in $G$
from $v_i$ to $v_j$, and similarly let $B_{j,i}$ be zero in $\R^{n_j\times l_i}$
if there is not an edge from $e_i$ to $v_j$. We then have
\[\bg(x_{k+1})\approx\bA\bg(x_k)+\bB\bh(u_k),\]
where
\begin{align*}
\bA&=\begin{bmatrix}
\bA_{1,1}&\bA_{1,2}&\cdots&\bA_{1,\nu}\\
\bA_{2,1}&\bA_{2,2}&\cdots&\bA_{2,\nu}\\
\vdots&\vdots&\ddots&\vdots\\
\bA_{\nu,1}&\bA_{\nu,2}&\cdots&\bA_{\nu,\nu}
\end{bmatrix},\\
\bB&=\begin{bmatrix}
\bB_{1,1}&\bB_{1,2}&\cdots&\bB_{1,\mu}\\
\bB_{2,1}&\bB_{2,2}&\cdots&\bB_{2,\mu}\\
\vdots&\vdots&\ddots&\vdots\\
\bB_{\nu,1}&\bB_{\nu,2}&\cdots&\bB_{\mu,\nu}
\end{bmatrix}.
\end{align*}

In the case where a reduced-order model is sought in the applications of the DMDc,
then for each $j=1,\ldots,\nu$, the process yields matrices
\begin{align*}
\tilde{\bA}_{j,j}=\hat{\bU}_j^*\bA_{j,j}\hat{\bU}_j,\\
\bar{\bA}_{j,k_j(i)}=\hat{\bU}_j^*\bA_{j,k_j(i)},\\
\tilde{\bB}_{j,\ell_j(i)}=\hat{\bU}_j^*\bB_{j,\ell_j(i)},
\end{align*}
where $\hat{U}_j$ is the matrix of left singular vectors in the SVD of
$\left[\begin{smallmatrix}\bZ_j\\\bGamma_j\end{smallmatrix}\right]$ that was
chosen in the DMDc application. These matrices satisfy
\begin{align*}
\tilde{\bg}_j(x_{k+1,j})\approx\text{}&\tilde{\bA}_{j,j}\tilde{\bg}_j(x_{k,j})+\sum_{i=1}^{\alpha_j}\bar{\bA}_{j,k_j(i)}\bg_{k_j(i)}(x_{k,k_j(i)})\\
&+\sum_{j=1}^{\beta_j}\tilde{\bB}_{j,\ell_j(i)}\bh_{\ell_j(i)}(u_{k,\ell_j(i)}),
\end{align*}
where $\tilde{\bg}_j:=\hat{\bU}_j^*\bg_j$. If we then let $\tilde{\bA}_{j,k}:=\bar{\bA}_{j,k}\hat{\bU}_k$, we have
\begin{align*}
\tilde{\bg}_j(x_{k+1,j})\approx\text{}&\tilde{\bA}_{j,j}\tilde{\bg}_j(x_{k,j})+\sum_{i=1}^{\alpha_j}\tilde{\bA}_{j,k_j(i)}\tilde{\bg}_{k_j(i)}(x_{k,k_j(i)})\\
&+\sum_{j=1}^{\beta_j}\tilde{\bB}_{j,\ell_j(i)}\bh_{\ell_j(i)}(u_{k,\ell_j(i)}).
\end{align*}
Also letting $\tilde{\bA}_{j,i}$ be the zero matrix of the appropriate size if there is not an edge in $G$ from $v_i$ to $v_j$, we obtain a reduced-order approximation of the entire system:
\[\tilde{\bg}(x_{k+1})\approx\tilde{\bA}\tilde{\bg}(x_k)+\tilde{\bB}\bh(u_k),\]
where
\begin{align*}
\tilde{\bg}=\begin{bmatrix}\tilde{\bg}_1\\\tilde{\bg}_2\\\vdots\\\tilde{\bg}_\nu\end{bmatrix},\quad
\bA&=\begin{bmatrix}
\tilde{\bA}_{1,1}&\tilde{\bA}_{1,2}&\cdots&\tilde{\bA}_{1,\nu}\\
\tilde{\bA}_{2,1}&\tilde{\bA}_{2,2}&\cdots&\tilde{\bA}_{2,\nu}\\
\vdots&\vdots&\ddots&\vdots\\
\tilde{\bA}_{\nu,1}&\tilde{\bA}_{\nu,2}&\cdots&\tilde{\bA}_{\nu,\nu}
\end{bmatrix},\\
\text{and}\quad\bB&=\begin{bmatrix}
\tilde{\bB}_{1,1}&\tilde{\bB}_{1,2}&\cdots&\tilde{\bB}_{1,\mu}\\
\tilde{\bB}_{2,1}&\tilde{\bB}_{2,2}&\cdots&\tilde{\bB}_{2,\mu}\\
\vdots&\vdots&\ddots&\vdots\\
\tilde{\bB}_{\nu,1}&\tilde{\bB}_{\nu,2}&\cdots&\tilde{\bB}_{\mu,\nu}
\end{bmatrix}.
\end{align*}

\section{Examples}\label{sec:examples}

We now illustrate the Network DMDc and its main advantages with a few examples.

\begin{example}[A simple linear example]\label{ex:simple1}
In this example, we illustrate in detail how the Network DMDc is applied to
the simple network described in Section \ref{ssec:networksystems} and depicted
in Figure \ref{fig:simplenetwork}. Assume all state and input spaces are
the real line $\R$, and that the dynamics is linear so that
\begin{align*}
T_{v_1}(x_{v_1},x_{v_2},x_{e_1})&=a_{1,1}x_{v_1}+a_{1,2}x_{v_2}+b_1x_{e_1}\text{, and}\\
T_{v_2}(x_{v_2},x_{e_2})&=a_{2,2}x_{v_2}+b_2x_{e_2}
\end{align*}
for fixed $a_{i,j},b_i\in\R$. Specifically, we set these constants as follows:
\[(a_{1,1},a_{1,2},a_{2,2},b_1,b_2)=(1.2,-0.5,0.8,1,1).\]
We shall see that the Network DMDc algorithm can recover the linear dynamics
of the system.

To perform Network DMDc, we simulate the system for $3$ time steps, yielding
data for $4$ total time instances. With
the notation of Section \ref{ssec:netdmdc}, we start with the initial state
$(x_{1,1},x_{1,2})=(2,5)$ and use the randomly generated sequences
of input values $(u_{1,1},u_{2,1},u_{3,1})=(0.2,0.4,0.8)$ (for $e_1$) and
$(u_{1,2},u_{2,2},u_{3,2})=(0.3,0.1,0.3)$ (for $e_2$).
Letting all observables $\bg_j,\bh_j$ be the identity, and $x_j=w_{j-1}$
for $j=2,3,4$, we get the following data matrices:
\begin{equation*}
\begin{aligned}[c]
\bZ_1&=[2\enspace 0.1\enspace {-1.63}]\\
\bY_1&=[0.1\enspace {-1.63} \enspace {-2.926}]\\
\bGamma_1&=\begin{bmatrix}5&4.3&3.54\\ 0.2&0.4&0.8\end{bmatrix}
\end{aligned}
\qquad
\begin{aligned}[c]
\bZ_2&=[5\enspace 4.3\enspace 3.54]\\
\bY_2&=[4.3\enspace 3.54\enspace 3.132]\\
\bGamma_2&=[0.3\enspace 0.1\enspace 0.3].
\end{aligned}
\end{equation*}
Performing DMDc on the triples $(\bZ_j,\bY_j,\bGamma_j)$ yields
\begin{align*}
[a_{1,1}\enspace a_{1,2}\enspace b_1]&=\bY_1\begin{bmatrix}\bZ_1\\\bGamma_1\end{bmatrix}^\dagger
=[1.2\enspace {-0.5}\enspace 1]\\
[a_{2,2}\enspace b_2]&=\bY_2\begin{bmatrix}\bZ_2\\\bGamma_2\end{bmatrix}^\dagger
=[0.8\enspace 1],
\end{align*}
thus recovering the dynamics.

Note that applying the regular DMDc algorithm to the above data does not
accurately recover the dynamics. However, regular DMDc can possibly
recover them if it is performed on data from a simulation spanning $5$
time instances. We explain why this happens after the next example, which
is a more dramatic illustration of this phenomenon.
\end{example}

\begin{example}[Circular networks]
We now consider a circular network with
state vertices $v_1,\ldots,v_n$ such that for $j=1,\ldots,n-1$, there is a
directed edge from $v_j$ to $v_{j+1}$, and there is another edge from
$v_n$ to $v_1$. Additionally, the network has input vertices $e_1,\ldots,e_m$,
$(m\leq n)$ each having a single edge pointing to one state vertex. On
the other hand, there is at most one edge pointing to a given state vertex
from an input vertex. Figure \ref{fig:circle} illustrates an example of this type of
network. As in the previous example, $P_w=\R$ for all $w\in N\cup I$, and the
dynamics is linear so that for a given $v_j\in N$, we have
\[T_{v_j}(x_{v_j},x_{v_{j-1(\text{mod\,}n)}},x_{e_k})
=a_jx_{v_j}+b_jx_{v_{j-1(\text{mod\,}n)}}+c_jx_{e_k}.\]
Here $k$ is the index such that $e_k$ is connected to $v_j$. If no input vertex
is connected to $v_j$, then the term $c_jx_{e_k}$ above is ignored.

\begin{figure}
\centering
\includegraphics[trim=40pt 30pt 40pt 30pt, clip, width=0.4\textwidth]{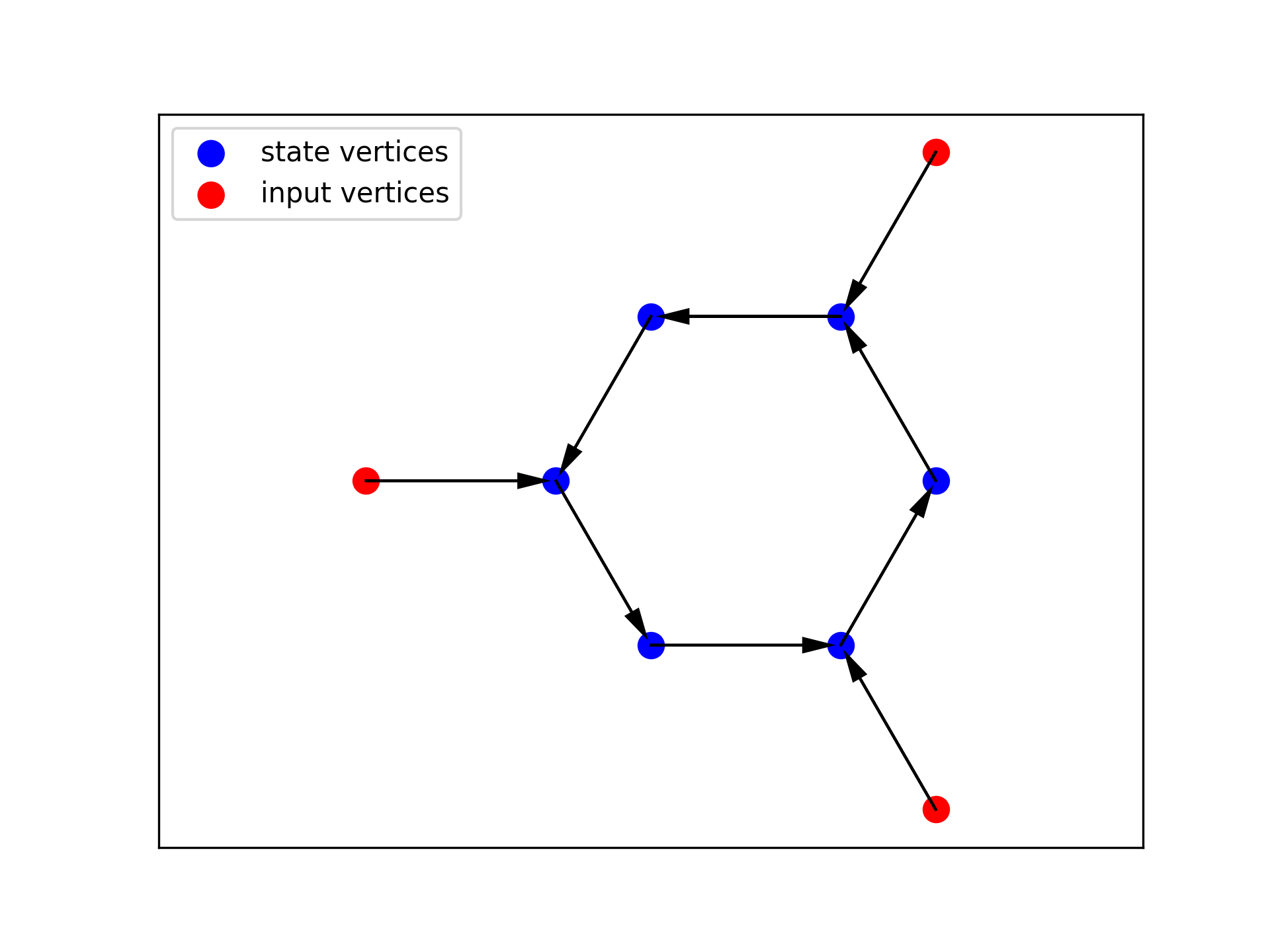}
\caption{Circular network with $6$ state vertices and $3$ input vertices}
\label{fig:circle}
\end{figure}

For several networks of this type having various sizes and parameters, we tested
the standard DMDc (or DMD in cases without inputs) and Network DMDc
algorithms by seeing how effectively both can recover their linear dynamics. In
particular, we performed simulations of the network of various time lengths,
resulting in triples
$(x_1,u_1,w_1),\ldots,(x_{m},u_{m},w_{m})\in\M\times\U\times\M$
with $x_j=w_{j-1}$ for $j\geq2$, and where the input values $u_k$
were randomly sampled from a uniform distribution on a finite interval. Then
for each simulation, we obtained the matrices $[\bA\enspace\bB]$ that result from applying
both algorithms to the data $(\bZ_j,\bY_j,\bGamma_j)$ following from letting
all observables $\bg_j,\bh_j$ be the identity. We then measured how close the matrices were to the
matrix giving the actual dynamics of the network. The measure of closeness used was the
Frobenious norm of the difference of the two matrices.

The Network DMDc was consistently more accurate than the standard DMDc throughout
these tests, and in particular required simulations of significantly shorter length to
recover the actual dynamics to within a negligible error. For instance, Figures \ref{fig:circularsim}
and \ref{fig:circularsimlog} depict the average error in the models generated by the network and standard DMDc algorithms for $20$ circular networks having $50$ state vertices, every other
of which is connected to an input vertex, analogous to the smaller network in Figure
\ref{fig:circle}. For each network, the parameters $a_j,b_j,c_j$ were chosen randomly from a uniform
distribution on the interval $[-1,1]$, and the values for the input vertices were
chosen from a uniform distribution on $[-10,10]$. (Note that for each simulation, the input data for the standard DMDc was the same as that for the network
DMDc, so that the accuracy comparison was fair.)

\begin{figure}
\centering
\includegraphics[trim=20pt 5pt 20pt 15pt, clip, width=0.4\textwidth]{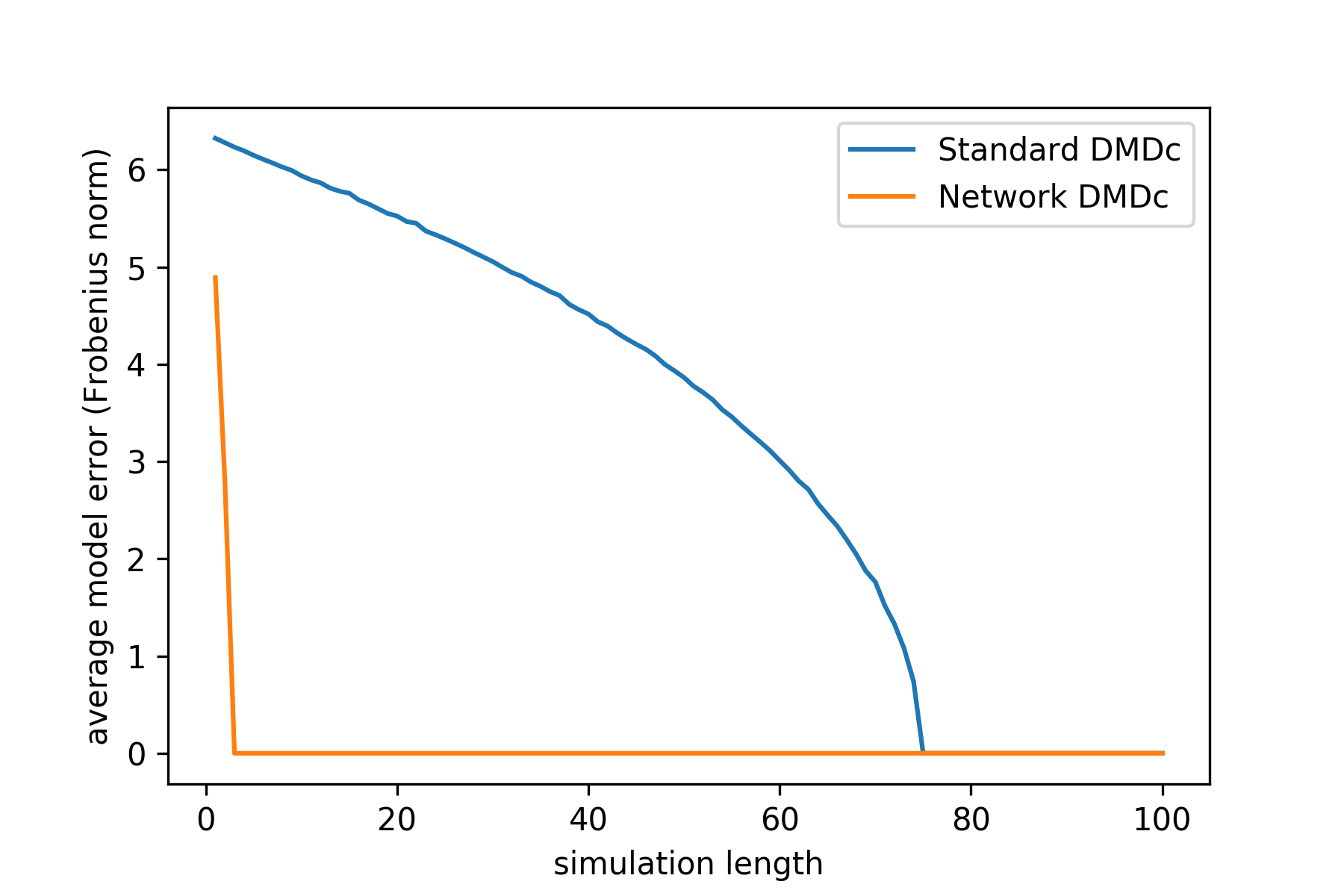}
\caption{The average error in the results of the standard and network DMDc applied to $20$ circular networks of $50$ state vertices $25$ input vertices}
\label{fig:circularsim}
\end{figure}

\begin{figure}
\centering
\includegraphics[trim=5pt 5pt 5pt 15pt, clip, width=0.4\textwidth]{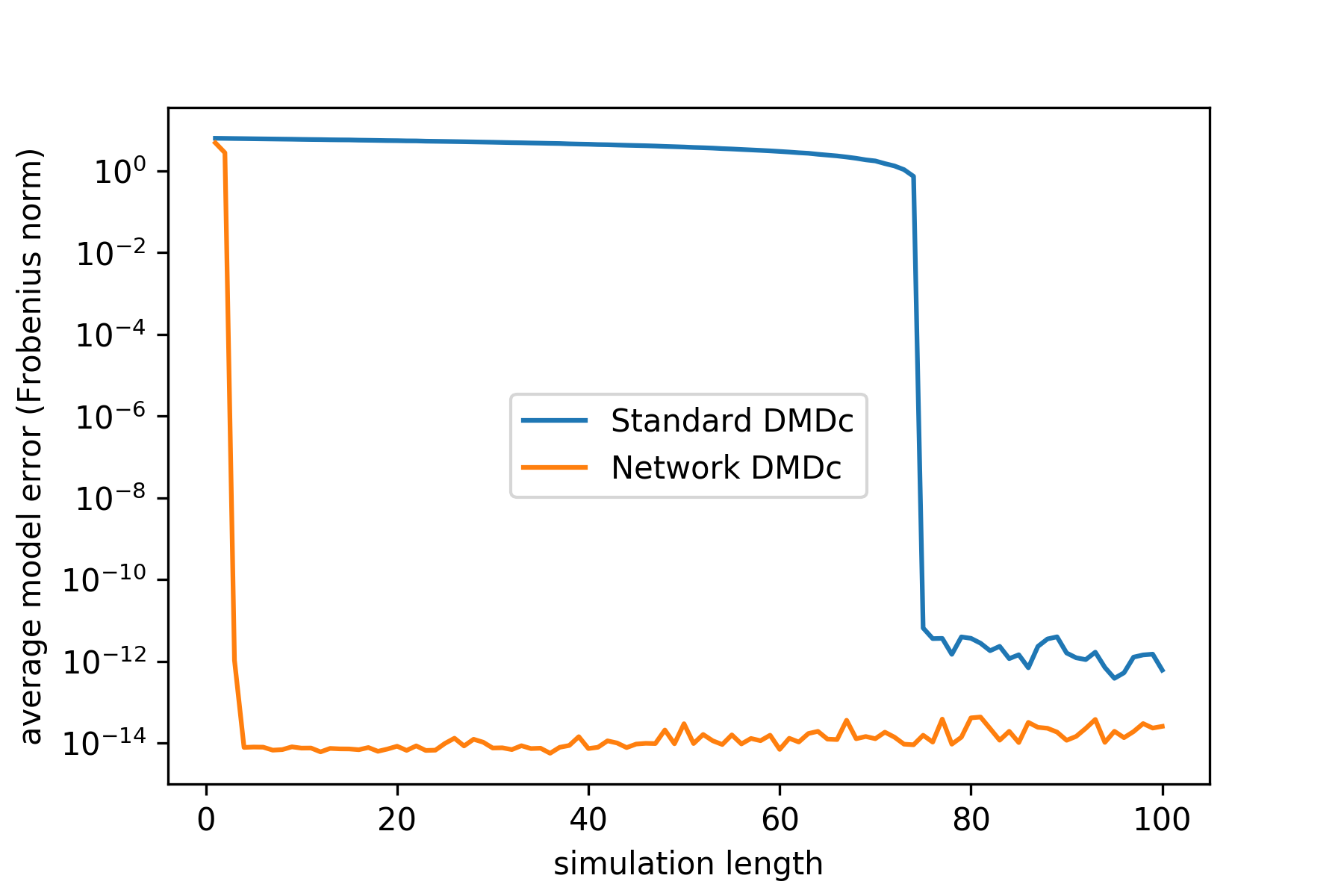}
\caption{Log scale of Figure \ref{fig:circularsim}}
\label{fig:circularsimlog}
\end{figure}

For the simulations depicted by Figures \ref{fig:circularsim} and
\ref{fig:circularsimlog}, the standard DMDc algorithm required a simulation
length of $m=75$ (that is, concluding with the calculation of $w_m=T(x_m,u_m)$)
to recapture the actual dynamics to within negligible error, whereas the
Network DMDc algorithm required only a simulation length of $m=3$. Note that $75$ is
the dimension of each system defined by one entire network, and $3$ is the
maximum dimension of each ``local subsystem'' consisting of one state vertex
$v$ and all the state and input vertices having edges pointing to $v$. 
(Recall that the Network DMDc is based on applying the regular DMDc to
each such local subsystem.) This makes sense because, in the first case, a
simulation length of $75$ would make each matrix
$\bOmega=\left[\begin{smallmatrix}\bZ\\\bGamma\end{smallmatrix}\right]$
a square matrix, and in the second case, a simulation length of $3$ would
make the matrices
$\left[\begin{smallmatrix}\bZ_j\\\bGamma_j\end{smallmatrix}\right]$
square or have more columns than rows. It is then possible that both
matrices have right inverses, in which case both algorithms would yield
the correct dynamics, up to computational error. (Problems could arise
if the data are correlated in a way so as to make the rows of $\bOmega$
or $\left[\begin{smallmatrix}\bZ_j\\\bGamma_j\end{smallmatrix}\right]$
linearly dependent.) Note that this also occurs for Example
\ref{ex:simple1}.
\end{example}

\begin{example}[Erd\H{o}s-Renyi random graphs]
To further establish the computational efficiency of Network DMDc,
we compared it to the standard DMD algorithm when applied to linear
systems whose network structures are Erd\H{o}s-Renyi random directed
graphs. That is, we considered networks produced according to the
following process: Fix a probability $p\in[0,1]$ and start with a
fixed number of vertices $n$. Then for every pair of vertices $(v,v')$,
include a directed edge from $v$ to $v'$ with probability $p$. These
choices of edges are determined independently. All vertices were assumed
to be state vertices. Also, $P_w=\R$ for all $w\in N$ and the scalars
determining the linear transition functions $T_w$ were chosen randomly
from a uniform distribution on $[-1,1]$.

The Network DMDc algorithm was performed as described in the previous
example on data generated by simulations of various time lengths of
these networks. The resulting matrix models were then compared to
those of the DMD (not DMDc since there are no control inputs). As in
the previous example, Network DMDc consistently outperformed DMD. For
example, for $20$ Erd\H{o}s-Renyi networks having $n=50$ vertices and
probability $p=0.05$ (an example is depicted in Figure
\ref{fig:erdosrenyinet}), the Network DMDc and DMD algorithms yielded
average model errors given in Figure \ref{fig:erdosrenyisim}. 

\begin{figure}
\centering
\includegraphics[trim=50pt 23pt 50pt 23pt, clip, width=0.4\textwidth]{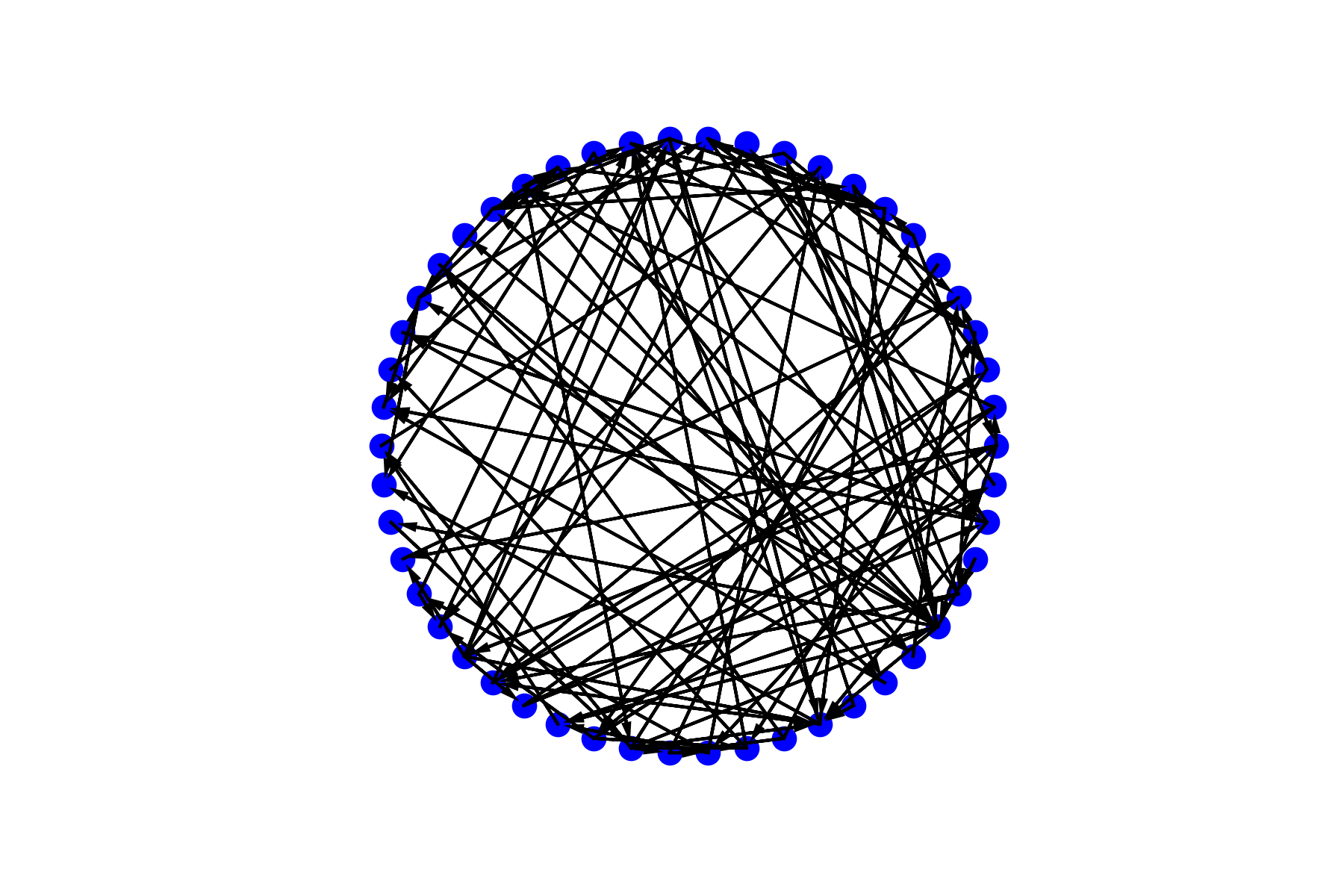}
\caption{Erd\H{o}s-Renyi random directed network with $n=50$ and $p=0.05$}
\label{fig:erdosrenyinet}
\end{figure}

\begin{figure}
\centering
\includegraphics[trim=6pt 5pt 6pt 9pt, clip, width=0.4\textwidth]{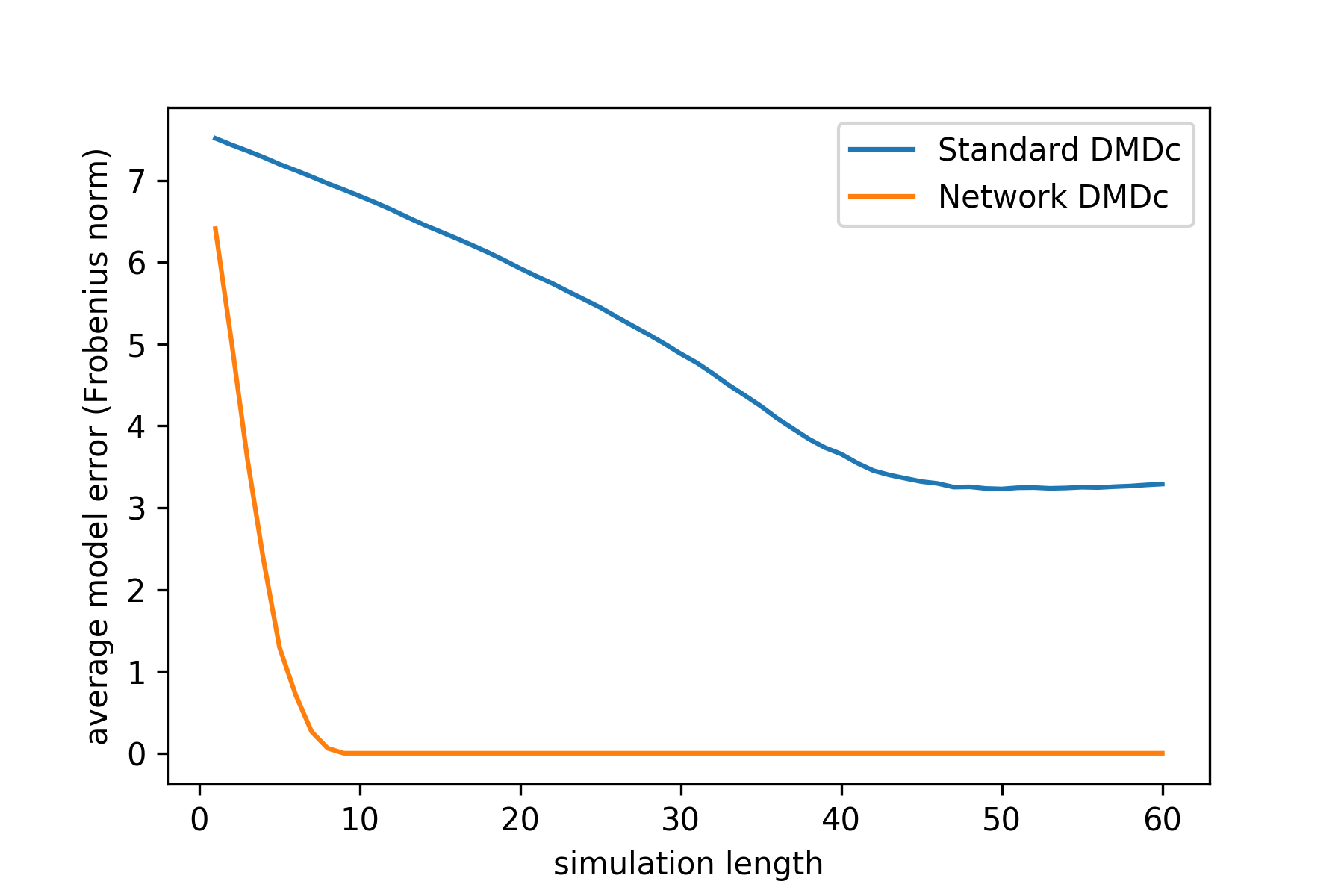}
\caption{Average error in the results of the DMD and Network DMDc applied to the $20$ Erd\H{o}s-Renyi random networks with $n=50$ and $p=0.05$}
\label{fig:erdosrenyisim}
\end{figure}

Once again, the Network DMDc algorithm recovered the linear dynamics with
data from simulations of relatively short length, namely the length matching
the dimension of each network's largest local subsystem. On the other hand, DMD was
unable to recover the dynamics accurately due to singular values too close to
$0$ (relative to the largest singular value) appearing in the data matrices
$\bOmega$ for large simulation lengths. This causes the pseudo-inverse, in our
case implemented by the \texttt{pinv()} function in the Python Numpy package
\cite{Scipy}, of the matrix difficult to numerically compute. Note that this is
liable to happen for Network DMDc as well, but is significantly mitigated if the
local subsystems of the network have low dimension.

\end{example}

The above examples demonstrate the computational benefit of Network
DMDc. Its central advantage is its exploitation of the network
structure in systems to decompose them into smaller subsystems, which
in turn lessens the burden and instability of the linear algebra
methods used in DMD and DMDc. We see in particular that these benefits
should generally manifest for network structures whose local subsystems
have a significantly smaller dimension than that of the system as a
whole (e.g., systems whose nodes have few incoming edges relative to the
size of the network).

\section{Conclusion}\label{sec:conclusion}

In this paper, we have developed a method of applying Koopman theory
and dynamic mode decomposition to networked control systems. In particular,
we have seen how to decompose the Koopman operator of a networked
system into lesser operators, on each of which we can apply the DMDc
algorithm. This allows us to obtain numerical approximations for the lesser
operators, which can then be composed to produce a linear model for the
entire system. We have seen through examples how this process can
improve the computation of the resulting models. By focusing on the dynamics
associated with each vertex separately and deliberately cutting out 
computation corresponding to dynamically unconnected components of the
system, Network DMDc can accurately recover linear dynamics of networked
systems with less data, and thus with a lower computational burden, than
standard DMD and DMDc. Additionally, Network DMDc lends itself naturally
to parallel computation, which can improve computational efficiency even
more. With its ability to also work with nonlinearities,
in addition to its possible use in tandem with distributed control algorithms,
Network DMDc has great potential to be used in the modeling and control
of complex interconnected systems.

\bibliographystyle{IEEEtran}
\bibliography{NetDMDcbib}

\end{document}